\documentclass[12pt]{a_t}
\usepackage{graphicx}
\usepackage{color}
\usepackage{url}

\usepackage{hyperref}
\hypersetup{
    colorlinks=true,
    linkcolor=black,
    citecolor=black,
    urlcolor=blue,
}
\usepackage{comment}
\usepackage{caption}
\usepackage{subcaption}

\usepackage{algorithm}
\usepackage{algorithmic}

\def\moh#1{{\color{black}#1}}

\begin{document}  

\year{2021}
\title{
О РЕШЕНИИ ВЫПУКЛЫХ MIN-MIN ЗАДАЧ С ГЛАДКОСТЬЮ И СИЛЬНОЙ ВЫПУКЛОСТЬЮ ПО ОДНОЙ ИЗ ГРУПП ПЕРЕМЕННЫХ И МАЛОЙ РАЗМЕРНОСТЬЮ ДРУГОЙ}%
\thanks{Работа выполнена  при поддержке Министерства науки и высшего образования Российской Федерации (госзадание) No. 075-00337-20-03, номер проекта
0714-2020-0005. Работа Е.Л. Гладина была также поддержана стипендией А.М. Райгородского в области численных методов оптимизации. Работа А.В. Гасникова была также частично поддержана грантом РФФИ 18-29-03071 мк.}

\authors{Е.Л.~ГЛАДИН\\
M.~АЛКУСА, канд.~физ.-мат.~наук\\
А.В.~ГАСНИКОВ, док.~физ.-мат.~наук\\
(Московский физико-технический институт, Долгопрудный)
}

\maketitle

\begin{abstract}
Статья посвящена некоторым подходам к решению выпуклых задач вида  min-min с гладкостью и сильной выпуклостью только по одной из двух групп переменных. Показано, что предложенные подходы, основанные на методе Вайды,
быстром градиентном методе и ускоренном градиентном методе с редукцией дисперсии, имеют линейную сходимость.
Для решения внешней задачи предлагается использовать методы Вайды, для решения внутренней (гладкой и сильно выпуклой)~--- быстрый градиентный метод. Ввиду важности для приложений в машинном обучении, отдельно рассмотрен случай, когда целевая функция является суммой большого числа функций. В этом случае вместо быстрого градиентного метода используется ускоренный градиентный метод с редукцией дисперсии. Приведены результаты численных экспериментов, иллюстрирующие преимущества предложенных процедур для задачи логистической регрессии, в которой есть априорное распределение на одну из двух групп переменных.

\end{abstract}

\section{Введение}
\moh{
Одним из основных направлений исследований численных методов выпуклой оптимизации в последнее десятилетие стало повсеместное распространение конструкции ускорения обычного градиентного метода, предложенной в 1983 г. Ю.Е. Нестеровым \cite{Nesterov1983}, на различные другие численные методы оптимизации. За последние 15 лет ускоренный метод был успешно перенесен на гладкие задачи условной выпуклой оптимизации, на задачи со структурой (в частности, так называемые композитные задачи), безградиентные  и рандомизированные методы (например, ускоренный градиентный метод с редукцией дисперсии для задач минимизации суммы функций \cite{Lan_book2020}). Также ускорение было успешно перенесено на методы, использующие старшие производные. Детали и более подробный обзор литературы можно найти в работе \cite{gasnikov2018book}.}


\moh{
Задачи оптимизации вида min-max и седловые  задачи широко изучались в литературе из-за их широкого спектра приложений в
статистике, машинном обучении, компьютерной графике, теории игр и других областях. В последнее время многие исследователи активно работают над темой ускоренных методов решения этих задач, учитывающих их структуру: \cite{Alkousa2020,Gladin2020,Jordan2020,Yuanhao2020,Zhongruo2020}~--- и это лишь некоторые из последних работ.  В некоторых приложениях существует задача, аналогичная задаче min-max, которая остается в значительной степени неизученной~--- это задача вида min-min:
\begin{equation}\label{main}
    \min_{x \in Q_x} \min_{y \in Q_y} F(x,y),
\end{equation}
где $Q_x\subset \mathbb{R}^d, Q_y\subset \mathbb{R}^n$~--- непустые компактные выпуклые  множества, размерность $d$ относительно небольшая ($d \ll n$), функция $F(x,y)$~--- выпуклая по совокупности переменных, а также $L$-гладкая и $\mu$-сильно выпуклая по $y$. Под $L$-гладкостью по $y$ понимается свойство
\begin{equation*}
    \|\nabla_y F(x, y)-\nabla_y F(x, y')\|_2 \leq L\|y-y'\|_2\quad \forall x \in Q_x, y,y' \in Q_y.
\end{equation*}
Такая постановка возникает, например, при поиске равновесий в транспортных сетях  \cite{Gasnikov20Models}.  В машинном обучении задачи такого типа соответствуют случаю, когда регуляризация применяется к одной из двух групп параметров модели (отсюда сильная выпуклость только по одной группе переменных из двух). Например, когда в датасете большая группа признаков являются разреженными, то регуляризация может использоваться только для весов модели, соответствующих этим признакам. В качестве ещё одного примера можно привести логистическую регрессию, в которой есть априорное распределение на часть параметров.
Задаче min-min посвящено несколько работ, среди которых \cite{Holderian2020,Jungers2011,Konur2017}. Например, в \cite{Holderian2020} авторы предложили новые алгоритмы для задач min-max, шаги которых настраиваются автоматически, но предложенные методы также применяются и к задачам  min-min.
}

\moh{
В данной статье мы рассматриваем два подхода к решению задачи \eqref{main}, имеющие линейную скорость сходимости.
Предлагается свести рассматриваемую задачу к совокупности вспомогательных задач (внутренней и внешней).
Внешняя задача (минимизация по $x$) решается методом Вайды (метод секущей плоскости) \cite{vaidya1989,vaidya1989new}.
}
\moh{
В случае, когда целевая функция $F$ простая, т.е. не является суммой большого количества функций, внутренняя задача (минимизация по $y$) решается быстрым градиентным методом для задач сильно выпуклой оптимизации.
В результате такого подхода приближённое решение задачи \eqref{main} может быть достигнуто за $\widetilde{\mathcal{O}}\left(d\right)$ вычислений $\partial_x F$ и $\widetilde{\mathcal{O}}\left(d\sqrt{\frac{L}{\mu}}\right)$ вычислений $\nabla_y F$, см. теорему \ref{theorem_approach_2}. Здесь и далее $\widetilde{\mathcal{O}}(\cdot) = \mathcal{O}(\cdot)$ с точностью до небольшой степени логарифмического множителя, обычно эта степень равна 1 или 2.
}
 
 \moh{
 Оптимизация
суммы большого количества функций в течение последних нескольких лет
является
предметом интенсивных исследований из-за широкого спектра приложений в машинном обучении, статистике, обработке изображений и других математических и инженерных приложениях. Поэтому мы отдельно рассматриваем случай, когда целевая функция $F$ представляет собой сумму (или среднее арифметическое) большого числа $m$ функций, в котором
использование быстрого градиентного метода для задач сильно выпуклой оптимизации потребовало бы вычисления градиентов $m$ слагаемых на каждом шаге, что может занимать много времени.
Вместо этого мы предлагаем использовать ускоренный градиентный метод с редукцией дисперсии \cite{Lan_book2020,Lan_NIPSpaper}, который также имеет линейную сходимость. В результате такого подхода, решение задачи  может быть достигнуто за $\widetilde{\mathcal{O}}\left(md\right)$ вычислений $\partial_x F$ и за $\widetilde{\mathcal{O}}\left(md + d \sqrt{\frac{mL}{\mu}}\right)$ вычислений $\nabla_y F$, см. теорему \ref{theorem_approach_3}.
 }
 
 \moh{
Используя два предложенных подхода, мы получаем линейную скорость сходимости для задачи 
min-min \eqref{main}. Отметим, что гладкость и сильная выпуклость требуется только по одной из двух групп переменных.
}

\moh{
Работа состоит из введения, заключения и 3 основных разделов. В разделе \ref{used_methods} мы приводим используемые алгоритмы и их сложность, а именно быстрый градиентный метод, метод Вайды (метод секущей плоскости)  и метод ускоренного градиентного спуска с редукцией дисперсии.
В разделе \ref{problem_stat_results} формулируется постановка задачи и приводятся
подходы к рассматриваемой задаче для различных случаев целевой функции, в одном из которых целевая функция является суммой или средним арифметическим большого числа функций. В разделе \ref{experiments} приводятся результаты вычислительных экспериментов и сравнение скорости работы предложенных подходов. Отметим, что полные доказательства теорем  \ref{delta_subgradient}, \ref{theorem_approach_2}, \ref{theorem_approach_3} и вспомогательного утверждения \ref{max_dot_product} приводятся в приложении к работе.
}

\section{Используемые алгоритмы}\label{used_methods}
\moh{
В этом разделе мы приводим алгоритмы, используемые в предлагаемых нами подходах к решению задачи \eqref{main}. Сначала приводится быстрый градиентный метод, затем метод Вайды (метод секущей плоскости) и, наконец, ускоренный градиентный метод с редукцией дисперсии.
}

\subsection{Быстрый градиентный метод}
\moh{
В  работе \cite{FGM_Gasnikov2019} предложен адаптивный алгоритм для решения следующей задачи оптимизации
\begin{equation}\label{problem_standard}
    f(y) \to \min_{y \in Q_y},
\end{equation}
где $Q_y \subset \mathbb{R}^n$~---непустое компактное выпуклое  множество, $f$~--- $L$-гладкая выпуклая функция.
Этот алгоритм, получивший название быстрого градиентного метода, позволяет ускорить сходимость обычного градиентного спуска с $\mathcal{O}\left(\frac{1}{N} \right)$ до $\mathcal{O}\left(\frac{1}{N^2} \right)$, где $N$~--- количество итерации алгоритма. Быстрый градиентный метод (не адаптивный вариант) приведён ниже как алгоритм \ref{alg:FGM}.
}

\begin{algorithm}
\caption{Быстрый градиентный метод \cite{FGM_Gasnikov2019}.
}
\label{alg:FGM}
\begin{algorithmic}[1]
	\REQUIRE Количество шагов $N$, начальная точка $y^0 \in Q_y$, параметр $L >0$.
	\STATE \textbf{0-шаг:} $z^{0}:=y^{0},\quad u^{0}:=y^{0},\quad \alpha_{0}:=0,\quad A_{0}:=0$.
	\FOR{$k = 0, 1, \ldots, N-1$}
	    \STATE Находим наибольший корень $\alpha_{k+1}$ такой, что
	    $A_{k}+ \alpha_{k+1} = L\alpha_{k+1}^2,$
		\STATE $A_{k+1} := A_k + \alpha_{k+1},$
		\STATE $z^{k+1} := \frac{\alpha_{k+1}u^k + A_k y^k}{A_{k+1}},$
		\STATE $u^{k+1} := \displaystyle{\arg\min_{y \in Q_y}}\left\{\alpha_{k+1} \left\langle \nabla f\left(z^{k+1}\right), y - z^{k+1}\right\rangle + \frac{1}{2}\|y - u^k\|_2^2  \right\}$,
		\STATE $y^{k+1} := \frac{\alpha_{k+1}u^{k+1} + A_k y^k}{A_{k+1}},$
	\ENDFOR
	\ENSURE $y^N.$
\end{algorithmic}
\end{algorithm}

\moh{
Следующая теорема дает оценку сложности (скорости сходимости) алгоритма \ref{alg:FGM}.
\begin{theorem}[{\cite{FGM_Gasnikov2019}}]\label{th:fgm}
Пусть функция $f: Q_y \to \mathbb{R}$ является $L$-гладкой и
выпуклой, тогда  алгоритм \ref{alg:FGM} возвращает такую точку $y^{N}$, что
$$
f\left(y^{N}\right)-f(y_*) \leq \frac{8LR^2}{(N+1)^2},
$$
где $y_*$~--- решение задачи \eqref{problem_standard}, $R^2 = \frac{1}{2}\|y^0 -y_*\|_2^2$.
\end{theorem}
}

\moh{Опишем далее технику рестартов (перезапусков) быстрого градиентного метода (алгоритм \ref{alg:FGM}), для случая $\mu$-сильно выпуклой функции.
}

\moh{Ввиду $\mu$-сильной выпуклости $f$ имеем
$$
\frac{\mu}{2}\|z-y\|_{2}^{2} \leq f(z)-\left(f(y)+  \langle\nabla f(y), z -y\rangle\right) \leq \frac{L}{2}\|z-y\|_{2}^{2}, \quad \forall y, z \in Q_y.
$$
Тогда после $N_1$ итераций алгоритма \ref{alg:FGM} с учётом теоремы \ref{th:fgm} получаем
\begin{equation}
    \frac{\mu}{2}\|y^{N_1} - y_* \|_2^2 \leq f\left(y^{N_1}\right)-f\left(y_*\right) \leq \frac{4L \|y^0 - y_* \|_2^2}{N_{1}^2},
\end{equation}
отсюда
$$
    \|y^{N_1} - y_*\|_2^2 \leq \frac{8L}{\mu N_1^2}\|y^0 - y_*\|_2^2.
$$
Поэтому, выбирая  $N_1 = \left\lceil 4 \sqrt{\frac{L}{\mu}}\right\rceil$, где $\left\lceil \cdot \right\rceil$~--- округление вверх, получим
$$
    \|y^{N_1} - y_*\|_2^2 \leq \frac{1}{2}\|y^0 - y_*\|_2^2.
$$
После этого выберем для алгоритма \ref{alg:FGM} в качестве точки старта $y^{N_1}$, и снова сделаем $N_1$ итераций, и т.д. Для достижения приемлемого качества решения можно выбрать количество рестартов алгоритма \ref{alg:FGM} (параметр $p$ алгоритма \ref{alg:restart}) следующим образом:
$$
    p = \left\lceil \frac{1}{2} \ln\left( \frac{\mu R^2}{\varepsilon}\right) \right\rceil.
$$
 В таком случае общее число итераций алгоритма \ref{alg:restart} будет
$$
    N = \left\lceil \frac{1}{2}\ln\left( \frac{\mu R^2}{\varepsilon}\right) \right\rceil \cdot \left\lceil 4 \sqrt{\frac{L}{\mu}} \right\rceil,
$$
т.е.
\begin{equation}\label{estimate_restat_alg}
    N = \mathcal{O}\left(\sqrt{\frac{L}{\mu}} \ln\left( \frac{\mu R^2}{\varepsilon}\right)  \right) = \widetilde{\mathcal{O}}\left( \sqrt{\frac{L}{\mu}}\right).
\end{equation}
}

\begin{algorithm}[htp]
	\caption{Быстрый градиентный метод для задач сильно выпуклой оптимизации, рестарты алгоритма \ref{alg:FGM}.}
	\begin{algorithmic}[1]
		\STATE {\bf Вход:} начальная точка $y^0 \in Q_y$, $L>0, \mu>0,$ точность решения $\varepsilon$, $R$, количество рестартов $p= \left\lceil \frac{1}{2}\ln\left( \frac{\mu R^2}{\varepsilon}\right) \right\rceil$.
		\FOR{$j = 1, \ldots, p$}
			\STATE Выполнить $N_j = \left\lceil 4 \sqrt{\frac{L}{\mu}} \right\rceil$ итераций алгоритма \ref{alg:FGM},
			\STATE $y^0 := y^{N_j}$.
		\ENDFOR
		\STATE {\bf Выход:} $\hat{y}:=y^{N_p}$.
	\end{algorithmic}
	\label{alg:restart}
\end{algorithm}

\subsection{Метод Вайды}
\moh{ 
Метод Вайды (метод секущей плоскости) был предложен Вайдой в \cite{vaidya1989,vaidya1989new} для решения следующей условной задачи оптимизации
\begin{equation}\label{problem_vaidya}
    f(x) \to \min_{x \in Q_x},
\end{equation}
где $Q_x \subset \mathbb{R}^d$~--- выпуклое компактное множество с непустой внутренностью, а целевая функция $f$, определённая на $Q_x$, непрерывна и выпукла.
}

\moh{
Пусть $P = \{x \in \mathbb{R}^d: \, Ax\geq b\}$~--- ограниченный $d$-мерный многогранник, где  $A \in \mathbb{R}^{m\times d}$ и $b \in \mathbb{R}^m$. Логарифмический барьер множества $P$ определяется как
$$
Barr(x) = -\sum_{i=1}^{m} \log \left(a_{i}^{\top} x-b_{i}\right),
$$
где  $a_{i}^{\top}$~--- $i$-я строка матрицы $A$. Гессиан  $H(x)$  функции $Barr(x)$ равен
$$
H(x) =\sum_{i=1}^{m} \frac{a_{i} a_{i}^{\top}}{\left(a_{i}^{\top} x-b_{i}\right)^{2}}
$$
Матрица $H(x)$ положительно определена для всех $x$ из внутренности $P$. Волюметрический барьер (volumetric barrier) $\mathcal{V}$ определяется как
$$
\mathcal{V}(x) = \frac{1}{2} \log \left(\operatorname{det}(H(x))\right),
$$
где $\operatorname{det}(H(x))$ обозначает детерминант $H(x)$. Будем называть точку минимума функции $\mathcal{V}$ на $P$ волюметрическим центром множества $P$.
}

\moh{Обозначим 
\begin{equation}\label{volumetric_barrier}
    \sigma_{i}(x)=\frac{a_{i}^{\top} \left(H(x)\right)^{-1} a_{i}}{\left(a_{i}^{\top} x-b_{i}\right)^{2}}, \quad 1 \leq i \leq m,
\end{equation}
тогда градиент волюметрического барьера $\mathcal{V}$ может быть записан как
$$
\nabla \mathcal{V}(x)=-\sum_{i=1}^{m} \sigma_{i}(x) \frac{a_i}{{a_i}^{\top}x - b_i}.
$$}
\moh{Пусть $\mathcal{Q}(x)$ определяется как
$$
\mathcal{Q}(x) = \sum_{i=1}^{m} \sigma_{i}(x) \frac{a_{i} a_{i}^{\top}}{\left(a_{i}^{\top} x-b_{i}\right)^{2}}.
$$
Заметим, что $\mathcal{Q}(x)$ положительно определена на внутренности $P$, а также $\mathcal{Q}(x)$ является хорошим приближением гессиана функции $\mathcal{V}(x)$, т.е. $\nabla^2 \mathcal{V}(x)$.
}

Метод Вайды производит последовательность пар $\left(A_t, b_t\right) \in \mathbb{R}^{m\times d}\times \mathbb{R}^{m}$ таких, что соответствующие многогранники содержат решение. В качестве начального многогранника, задаваемого парой $\left(A_0, b_0\right)$, обычно берётся симплекс (алгоритм может начинать с любого выпуклого ограниченного $n$-мерного многогранника, для которого легко вычислить волюметрический центр~--- например, с $n$-прямоугольника).

Фиксируем небольшую константу $\gamma \geq 0.006$. Пусть $x_k$ ($k\geq 0$) обозначает волюметрический центр многогранника, заданного парой $\left(A_k, b_k\right)$, и пусть для него вычислены величины $\left\{\sigma_i(x_k)\right\}_{1\leq i \leq m}$ (см. \eqref{volumetric_barrier}). Следующий многогранник $\left(A_{k+1}, b_{k+1}\right)$ получается из текущего в результате либо присоединения, либо удаления ограничения:
\begin{enumerate}
    \item Если для некоторого $i \in \{1,\ldots, m\}$ выполняется $\sigma_i(x_k) = \min\limits_{1\leq j\leq m}\sigma_j(x_t) < \gamma$, тогда $\left(A_{k+1}, b_{k+1}\right)$ получается исключением $i$-й строки из $\left(A_k, b_k\right)$.
    \item Иначе (если $\min\limits_{1\leq j\leq m}\sigma_j(x_k) \geq \gamma$) оракул, вызванный в текущей точке $x_k$, возвращает вектор $c_k$ такой, что $f(x) \leq f(x_k)\ \forall x \in \left\{z \in Q_x: c_k^{\top} z\geq c_k^{\top} x_k \right\}$, т.е. $c_k \in -\partial f(x_k)$. Выберем $\beta_k \in \mathbb{R}$ таким, что
    $$
    \frac{c_k^{\top} \left(H(x_k)\right)^{-1} c_k}{\left(x_{k}^{\top} c_k-\beta_{k}\right)^{2}}=\frac{1}{5} \sqrt{\gamma}.
    $$
    Определим $\left(A_{k+1}, b_{k+1}\right)$ добавлением строки $\left(c_k, \beta_k\right)$ к $\left(A_k, b_k\right)$.
\end{enumerate}
Волюметрический барьер $\mathcal{V}_k$ является самосогласованной функцией, поэтому может быть эффективно минимизирован методом Ньютона. Достаточно одного шага метода Ньютона для $\mathcal{V}_k$, сделанного из $x_{k-1}$. Подробности и анализ метода Вайды можно найти в \cite{vaidya1989,vaidya1989new} и \cite{Bubeck}.

\moh{
Следующая теорема дает оценку сложности алгоритма Вайды.
}
\begin{theorem}\label{th:vaidya}
    Пусть выпуклое компактное множество $Q_x \subseteq \mathcal{B}_{\mathcal{R}}$ и $\mathcal{B}_{\rho} \subseteq \left\{x \in Q_x: f(x)-f\left(x_* \right) \leqslant \varepsilon \right\}$, где $\mathcal{B}_{\rho}$~--- евклидов шар радиуса $\rho$. Тогда метод Вайды находит $\varepsilon$-решение задачи \eqref{problem_vaidya} за $\mathcal{O} \left( d \log \frac{d \mathcal{R}}{\rho} \right)$ шагов.
\end{theorem}

\begin{remark}\label{mat_inversion}
    Помимо вычисления субградиента, в стоимость итерации метода Вайды входит стоимость обращения матрицы размера $d \times d$ и решения системы линейных уравнений.
\end{remark}

\subsection{Ускоренный градиентный метод с редукцией дисперсии}
\moh{
Рассмотрим задачу
\begin{equation}\label{constrained_problem}
    f(y) \to \min_{y \in Q_y},
\end{equation}
где $Q_y \subseteq \mathbb{R}^n$~--- замкнутое выпуклое множество, а целевая функция $f$ представляет собой сумму (или среднее арифметическое) большого числа $m$ гладких 
выпуклых функций $f_i$, т.е. $f(y) = \frac{1}{m}\sum\limits_{i=1}^{m}f_i(y)$. При решении \eqref{constrained_problem} с помощью быстрого градиентного метода для задач сильно выпуклой оптимизации (алгоритм \ref{alg:restart}) потребуется вычислять градиент $m$ функций на каждой итерации, что очень дорого. Поэтому предпочтительнее вместо алгоритма \ref{alg:restart} использовать рандомизированный градиентный метод, а именно ускоренный градиентный метод с редукцией дисперсии, также называемый Varag \cite{Lan_book2020,Lan_NIPSpaper}.
Приведённый ниже алгоритм \ref{varag} представляет собой ускоренный градиентный метод с редукцией дисперсии (Varag) для гладкой сильно выпуклой задачи оптимизации конечной суммы \eqref{constrained_problem}. Этот алгоритм был предложен Г. Ланом и др. в \cite{Lan_NIPSpaper}.
}
\moh{
Мы предполагаем, что для каждого $i \in \{1,\ldots,m\}$, существует $L_i>0$ такое, что
$$
\|\nabla f_i(y) - \nabla f_i(z)\|_2 \leq L_i \|y-z\|_2, \quad \forall y,z \in Q_y.
$$
Ясно, что $f$ имеет липшицев градиент с константой не более $L := \frac{1}{m}\sum\limits_{i=1}^{m}L_i$. Мы также предполагаем, что целевая функция $f$ сильно выпуклая с константой $\mu \geq 0$, т.е.
$$
f(z) \geq f(y)+\langle\nabla f(y), z-y\rangle + \frac{\mu}{2}\|y-z\|_2, \quad \forall y, z \in Q_y.
$$
}

\moh{
Алгоритм Varag содержит вложенные циклы~--- внешний и внутренний (индексируемые переменными $s$ и $t$, соответственно).
На каждой итерации внешнего цикла вычисляется полный градиент $\nabla f(\tilde{y}) $ в точке $\tilde{y}$, который затем используется во внутреннем цикле для определения оценок градиента $G_t$.
Каждая итерация внутреннего цикла требует информацию о градиенте только одного случайно выбранного слагаемого $f_{i_t}$ и содержит три основные последовательности: $\{\underline{y}_t\}, \{y_t\}$ и $\{\bar{y}_t\}$.} 

Обозначим $s_0 := \lfloor\log_2 m\rfloor + 1$,  где $\lfloor \cdot \rfloor$~--- округление вниз. Параметры алгоритма \ref{varag} $\{q_1, \ldots, q_m\}$, $\{\theta_t\}$, $\{\alpha_s\}$, $\{\gamma_s\}$, $\{p_s\}$ и $\{T_s\}$  описываются следующим образом:
\begin{itemize}
    \item Вероятности $q_i= \frac{1}{\sum\limits_{i=1}^{m}L_i}L_i, \; \forall i \in \{1, \ldots, m\}$.
    \item Веса $\{\theta_t\}$ при $1 \leq s \leq s_{0}$ или $s_0 < s \leq s_0  + \sqrt{\frac{12 L}{m \mu}}-4, m<\frac{3 L}{4 \mu}$ равны
    \begin{equation}\label{weights_theta_1}
        \theta_{t}=\left\{\begin{array}{ll}
        \frac{\gamma_{s}}{\alpha_{s}}\left(\alpha_{s}+p_{s}\right), & 1 \leq t \leq T_{s}-1, \\
        \gamma_s/\alpha_s, & t=T_s.
        \end{array}\right.
    \end{equation}
    В остальных случаях они равны
    \begin{equation}\label{weights_theta_2}
        \theta_{t}=\left\{
        \begin{array}{ll}
            \Gamma_{t-1}-\left(1-\alpha_{s}-p_{s}\right) \Gamma_{t}, & 1 \leq t \leq T_{s}-1, \\
            \Gamma_{t-1}, & t=T_{s},
        \end{array}\right.
    \end{equation}
    где $\Gamma_{t}=\left(1+\mu \gamma_{s}\right)^{t}$.
    \item Параметры $\{T_s\}, \{\gamma_s\}$ и $\{p_s\}$ определяются как
    \begin{equation}\label{gamma_p_T}
        T_s=\left\{\begin{array}{ll}
        2^{s-1}, & s \leq s_0 \\
        T_{s_0}, & s > s_0
        \end{array}, \quad  \gamma_s=\frac{1}{3 L \alpha_s},\right.  \quad \text {и} \;\; p_s=\frac{1}{2}.
    \end{equation}
    \item Наконец,
    \begin{equation}\label{alpfa}
        \alpha_{s}=\left\{\begin{array}{ll}
        \frac{1}{2}, & s \leq s_{0}, \\
        \max \left\{\frac{2}{s-s_{0}+4}, \min \left\{\sqrt{\frac{m \mu}{3 L}}, \frac{1}{2}\right\}\right\}, & s>s_{0}.
        \end{array}\right.
    \end{equation}
\end{itemize}

\moh{\begin{algorithm}[htp]
\caption{
Ускоренный градиентный метод с редукцией дисперсии (Varag) \cite{Lan_NIPSpaper}.}
\label{varag}
\begin{algorithmic}[1]
	\REQUIRE  $y^0 \in Q_y, \{T_s\}, \{\gamma_s\}, \{\alpha_s\}, \{p_s\}, \{\theta_t \}$ и распределение вероятностей  $\{q_1, \ldots, q_m\}$ на $\{1, \ldots, m\}$.
	\STATE $\tilde{y}^0:=y^0$. 
	\FOR{$s = 1, 2, \ldots, $}
	    \STATE $\tilde{y}:=\tilde{y}^{s-1}, \tilde{g}:=\nabla f(\tilde{y})$. 
	    \STATE $y_0:=y^{s-1}, \bar{y}_0=\tilde{y}, T: = T_s$.
	    \FOR{$t = 1, 2, \ldots, T$}
	    \STATE Выбрать $i_t \in \{1, \ldots, m\}$ случайным образом согласно $\{q_1, \ldots, q_m\}$.
	    \STATE $\underline{y}_t:= \frac{1}{\left(1+\mu \gamma_s(1-\alpha_s)\right)} \left[(1+\mu \gamma_s)(1-\alpha_s - p_s) \bar{y}_{t-1}+\alpha_{s} y_{t-1}+\left(1+\mu \gamma_{s}\right) p_{s} \tilde{y}\right]$.
		\STATE $G_t:=\frac{1}{(q_{i_t} m)}\left(\nabla f_{i_t}\left(\underline{y}_{t}\right)-\nabla f_{i_t}(\tilde{y})\right)+\tilde{g}$.
		\STATE $y_t:=\arg \min_{y \in Q_y}\left\{\gamma_s\left(\left\langle G_t, y\right\rangle+\frac{\mu}{2} \|\underline{y}_{t} -y\|_2^2\right)+\frac{1}{2}\|y_{t-1}- y\|_2^2\right\}$.
		\STATE $\bar{y}_t:=\left(1-\alpha_s - p_s\right) \bar{y}_{t-1}+\alpha_s y_t + p_s \tilde{y}$.
	    \ENDFOR
	    \STATE $y^s := y_T, \tilde{y}^{s}:= \frac{1}{\sum_{t=1}^{T} \theta_t}\sum_{t=1}^{T}\left(\theta_{t} \bar{y}_{t}\right)$.
	    \ENDFOR
\end{algorithmic}
\end{algorithm}
}

\moh{
Следующий результат дает оценку сложности алгоритма \ref{varag}.
\begin{theorem}[\cite{Lan_NIPSpaper}]\label{theorem_estimate_varag}
Если параметры алгоритма \ref{varag} $ \{\theta_t\}, \{\alpha_s\}, \{\gamma_s\}, \{p_s\}$ и $\{T_s\}$, заданы согласно формулам \eqref{weights_theta_1}, \eqref{weights_theta_2}, \eqref{gamma_p_T} и \eqref{alpfa}, то общее количество вычислений градиентов функций $f_i$, выполняемых алгоритмом \ref{varag} для нахождения стохастического $\varepsilon$-решения задачи \eqref{constrained_problem}, ограничено
\begin{equation}\label{estimate_varag}
    \widetilde{\mathcal{O}}\left(m +\sqrt{\frac{m L}{\mu}}\right),
\end{equation}
где в \eqref{estimate_varag} $\widetilde{\mathcal{O}}() = \mathcal{O}()$ с точностью до логарифмического множителя по $m, L, \mu, \varepsilon$ и $D_0 =  2\left(f(y^0) - f(y_*)\right) + \frac{3L}{2}\|y^0 - y_*\|_2^2$, где $y_*$~--- решение задачи \eqref{constrained_problem}.
\end{theorem}
}

\moh{
Из теоремы \ref{theorem_estimate_varag} видно, что алгоритм \ref{varag} достигает хорошо известной оптимальной линейной скорости сходимости.
}

\section{Постановка задачи и полученные результаты}\label{problem_stat_results}
Рассмотрим задачу
\begin{equation}\label{main4}
    \min_{x \in Q_x} \min_{y \in Q_y} F(x,y),
\end{equation}
где $Q_x \subset \mathbb{R}^d,\ Q_y \subset \mathbb{R}^n$~--- непустые компактные выпуклые множества, размерность $d$ относительно небольшая ($d \ll n$), функция $F(x,y)$~--- выпуклая по совокупности переменных, а также $L$-гладкая и $\mu$-сильно выпуклая по $y$. Под $L$-гладкостью по $y$ понимается свойство
\begin{equation*}
    \|\nabla_y F(x, y)-\nabla_y F(x, y')\|_2 \leq L\|y-y'\|_2\quad \forall x \in Q_x, y,y' \in Q_y.
\end{equation*}

Введём функцию
\begin{equation}\label{aux2}
    f(x) = \min_{y \in Q_y} F(x,y).
\end{equation}

Задачу \eqref{main4} можно переписать следующим образом:
\begin{equation}\label{min_g}
    f(x) \to \min_{x \in Q_x}
\end{equation}
При решении \eqref{min_g} некоторым итерационным методом необходимо на каждом его шаге решать вспомогательную задачу \eqref{aux2}, чтобы приближённо находить субградиент $\partial f(x)$. Обратимся к следующему определению.
\begin{definition}[(\cite{polyak1983intro}, с. 123)]
Пусть $\delta \geq 0,\ Q_x \subseteq \mathbb{R}^d$~--- выпуклое множество, $f: Q_x \to \mathbb{R}$~--- выпуклая функция. Вектор $g \in \mathbb{R}^d$ называется $\delta$-субградиентом $f$ в точке $x' \in Q_x$, если
\begin{equation*}
    f(x) \geq f(x') + \langle g, x-x' \rangle - \delta\quad \forall x \in Q_x.
\end{equation*}
Множество $\delta$-субградиентов $f$ в точке $x'$ обозначается $\partial_\delta f(x')$.
\end{definition}

Обозначим $\displaystyle D:= \max_{y, z \in Q_y} \|y-z\|_2,\ y(x):= \arg\min_{y \in Q_y} F(x, y)$. Следующая теорема говорит о том, как вычислить $\delta$-субградиент функции $f(x)$, приближённо решая вспомогательную задачу \eqref{min_g}.
\begin{theorem}\label{delta_subgradient}
    Пусть найден такой $\tilde{y} \in Q_x$, что $F(x, \tilde{y}) - f(x) \leq \varepsilon$, тогда
    \begin{equation*}
        \partial_{x} F(x, \tilde{y}) \in \partial_\delta f(x),\quad \delta = \left(LD+\left\| \nabla_y F\left(x, y(x)\right) \right\|_2^2\right) \sqrt{\frac{2 \varepsilon}{\mu}}.
    \end{equation*}
\end{theorem}
Эта теорема непосредственно следует из двух утверждений:
\begin{statement}\label{max_dot_product}
    Пусть $g: Q_y \to \mathbb{R}$~--- $L$-гладкая $\mu$-сильно выпуклая функция, точка $\tilde{y} \in Q_y$ такова, что $g(\tilde{y}) - g(y_*) \leq \varepsilon$, тогда
    \begin{equation*}
        \max_{y \in Q_y}\left\langle \nabla g(\tilde{y}), \tilde{y}-y\right\rangle \leq \delta,\quad \delta = \left(LD+\left\| \nabla g\left(y_*\right) \right\|_2^2\right) \sqrt{\frac{2 \varepsilon}{\mu}},
    \end{equation*}
    где $y_* = \arg\min_{y \in Q_y} g(y)$.
\end{statement}
\begin{statement}[(\cite{Gas2015Equi}, с. 12)]\label{statement_gas}
    Пусть найден такой $\tilde{y} \in Q_y$, что
    \begin{equation*}
        \max_{y \in Q_y}\left\langle \nabla_y F(x, \tilde{y}), \tilde{y}-y\right\rangle \leq \delta,
    \end{equation*}
    тогда $\partial_{x} F(x, \tilde{y}) \in \partial_\delta f(x)$.
\end{statement}
Интуитивно теорема \ref{delta_subgradient} говорит о том, что, решив вспомогательную задачу \eqref{aux2} достаточно точно, мы получим хорошее приближение субградиента $\partial f(x)$, которое может быть использовано для решения внешней задачи \eqref{min_g}. На этой идее основан предлагаемый подход к решению \eqref{main4}.

\begin{approach}[(основной случай)]\label{approach_1}
Внешняя задача \eqref{min_g} решается методом Вайды.  
Вспомогательная задача \eqref{aux2} решается быстрым градиентным методом для задач сильно выпуклой оптимизации (алгоритм \ref{alg:restart}).
\end{approach}

\begin{theorem}\label{theorem_approach_2}
Подход \ref{approach_1} позволяет получить $\varepsilon$-решение задачи \eqref{main4} после $\widetilde{\mathcal{O}} \left( d \right)$ вычислений $\partial_x F$ и обращений матриц размера $d \times d$, а также $\widetilde{\mathcal{O}} \left(d \sqrt{\frac{L}{\mu}} \right)$ вычислений $\nabla_{y} F$.
\end{theorem}

\begin{remark}
Обращение матриц появляется в сложности предлагаемого подхода из-за того, что оно производится на каждом шаге метода Вайды.
\end{remark}

\subsection{Минимизация суммы большого числа функций}
Пусть в условиях предыдущего раздела
\begin{equation}
    F(x, y) = \frac{1}{m} \sum_{i=1}^m F_i (x, y),
\end{equation}
где функции $F_i$ являются выпуклыми по совокупности переменных и $L_i$-гладкими по $y$, а $F$ является $\mu$-сильно выпуклой по $y$. Из этого следует, что $F$ является выпуклой по совокупности переменных и гладкой по $y$ с константой гладкости не более $L := \frac{1}{m}\sum\limits_{i=1}^{m}L_i$.

\begin{approach}[(сумма функций)]\label{approach_2}
Внешняя задача \eqref{min_g} решается методом Вайды. 
 Вспомогательная задача \eqref{aux2} решается ускоренным градиентным методом с редукцией дисперсии (алгоритм \ref{varag}). 
\end{approach}

\begin{theorem}\label{theorem_approach_3}
Подход \ref{approach_2} позволяет получить $\varepsilon$-решение задачи \eqref{main4} за $\widetilde{\mathcal{O}}\left(md\right)$ вычислений $\partial_x F_i$, $\widetilde{\mathcal{O}}\left(d\right)$ обращений матриц размера $d \times d$ и $\widetilde{\mathcal{O}}\left(dm + d\sqrt{\frac{m L}{\mu}} \right)$ вычислений $\nabla_y F_i$.
\end{theorem}

\section{Эксперименты}\label{experiments}
Рассмотрим модель логистической регрессии для задачи классификации. Ошибка модели с параметрами $w$ на обучающем объекте с вектором признаков $z$, принадлежащем классу $t \in \{-1, 1\}$, записывается как
\begin{equation*}
    \ell_z(w) = \log \left(1+e^{-t \langle w, z \rangle}\right).
\end{equation*}
Пусть параметры модели состоят из двух групп: $w=(x, y), x \in \mathbb{R}^d, y \in \mathbb{R}^n$, при чём на группу $y$ наложено гауссовское априорное распределение: $y \sim \mathcal{N}\left(0, \sigma^2 I_n \right)$, где $I_n$~--- единичная матрица размера $n$. Максимизация апостериорной вероятности приведёт к задаче
\begin{equation}\label{logreg}
    \min_{x \in Q_x} \min_{y \in Q_y} \left\{ F(x,y):= \frac{1}{m} \sum_{i=1}^m \ell_{z_i}(x,y) + \frac{1}{\sigma^2} \| y \|_2^2 \right\},
\end{equation}
где в качестве $Q_x$ и $Q_y$ можно взять евклидовы шары достаточно большого радиуса.

Будем решать задачу \eqref{logreg} при помощи подхода \ref{approach_2} и сравним его работу с работой метода Varag (алгоритм \ref{varag}). Заметим, что эта задача не является сильно выпуклой по совокупности переменных. Для такой постановки можно использовать Varag, задавая параметры $\theta_t$ по формуле \eqref{weights_theta_1}, а все остальные параметры по формулам для сильно выпуклого случая, положив $\mu = 0$, см. \cite{Lan_NIPSpaper}. При этом стохастическое $\varepsilon$-решение будет найдено за $\mathcal{O}\left( \sqrt{\frac{m D_0}{\varepsilon}} + m \log m \right)$ вычислений градиентов функций $F_i$, где $D_0 =  2\left(F(x^0, y^0) - F(x_*, y_*)\right) + \frac{3L}{2}\|(x^0, y^0) - (x_*, y_*)\|_2^2$, $(x_*, y_*)$~--- решение задачи \eqref{logreg}. Эта сублинейная оценка уступает предлагаемому нами подходу, см. теорему \ref{approach_2}.

Для экспериментов использовался датасет \href{https://www.csie.ntu.edu.tw/~cjlin/libsvmtools/datasets/binary.html#madelon}{madelon}, представленный 2000 объектами, имеющими 500 признаков. Был выбран небольшой коэффициент регуляризации $\frac{1}{\sigma^2} = 0.005$ и проведены эксперименты для двух размерностей $d$, равных 20 и 30.

\begin{figure}
     \centering
     \begin{subfigure}[b]{0.49\textwidth}
         \centering
         \includegraphics[width=0.99\textwidth]{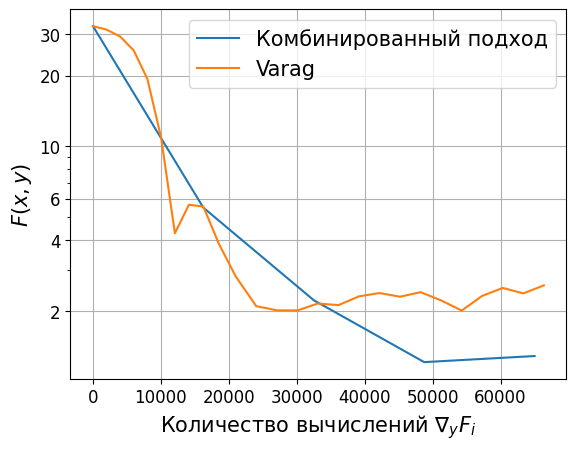}
         \caption{$d=20$}
         \label{fig:exp_d20}
     \end{subfigure}
     \hfill
     \begin{subfigure}[b]{0.49\textwidth}
         \centering
         \includegraphics[width=0.99\textwidth]{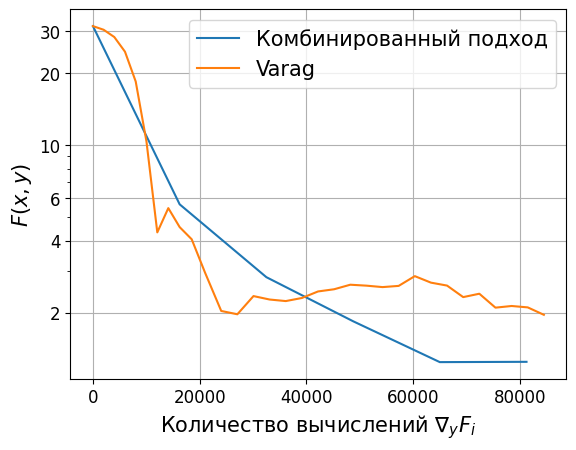}
         \caption{$d=30$}
         \label{fig:exp_d30}
     \end{subfigure}
        \caption{Результаты экспериментов.}
        \label{fig:exp}
\end{figure}

На рисунке \ref{fig:exp} отражены результаты эксперимента. По оси x откладывается количество вычислений градиентов $\nabla_y F_i$, которое для Varag совпадает с количеством вычислений $\nabla_x F_i$. Отметим, что наш подход требует меньше вычислений $\nabla_x F_i$, поскольку они выполняются только во внешнем цикле. Так, график \ref{fig:exp_d20} соответствует 4 итерациям внешнего цикла (т.е. 8000 вычислений $\nabla_x F_i$), а график \ref{fig:exp_d30}~--- 5 итерациям (т.е. 10000 вычислений $\nabla_x F_i$). В данном эксперименте подход \ref{approach_2} позволил достичь меньших значений целевой функции.

Исходный код и результаты экспериментов могут быть найдены в репозитории \url{https://github.com/egorgladin/min_min}.

\section{Заключение}
В статье рассмотрена задача вида min-min:
\begin{equation}\label{main__}
    \min_{x \in Q_x} \min_{y \in Q_y} F(x,y),
\end{equation}
где $Q_x\subset \mathbb{R}^d, Q_y\subset \mathbb{R}^n$~--- непустые компактные выпуклые  множества, размерность $d$ относительно небольшая ($d \ll n$), функция $F(x,y)$~--- выпуклая по совокупности переменных, а также $L$-гладкая и $\mu$-сильно выпуклая по $y$.

Предложено два подхода к решению задачи \eqref{main__}, в которых она сводится к совокупности вспомогательных задач (внутренней и внешней).
Внешняя задача (минимизация по $x$) решается методом Вайды, а внутренняя (минимизация по $y$)~--- быстрым градиентным методом для задач сильно выпуклой оптимизации или, если минимизируется сумма большого количества функций, ускоренным градиентным методом с редукцией дисперсии. Это позволяет достигать приближённого решения задачи \eqref{main__} за $\widetilde{\mathcal{O}}\left(d\right)$ вычислений $\partial_x F$ и $\widetilde{\mathcal{O}}\left(d\sqrt{\frac{L}{\mu}}\right)$ вычислений $\nabla_y F$, см. теорему \ref{theorem_approach_2}. Для сравнения, если бы задача \eqref{main__} была гладкой по совокупности переменных, то её решение при использовании только быстрого градиентного метода имело бы сложность $\mathcal{O}\left( \sqrt{\frac{LR^2}{\varepsilon}} \right)$, где $R$~--- расстояние от начального приближения до решения.
В случае суммы с $m$ слагаемыми, решение задачи  может быть достигнуто за $\widetilde{\mathcal{O}}\left(md\right)$ вычислений $\partial_x F$ и за $\widetilde{\mathcal{O}}\left(md + d \sqrt{\frac{mL}{\mu}}\right)$ вычислений $\nabla_y F$, см. теорему \ref{theorem_approach_3}.

Проведён численный эксперимент, в котором один из предлагаемых подходов применён к задаче логистической регрессии с регуляризацией, применяемой к одной из двух групп параметров модели. По сравнению с алгоритмом Varag, наш подход достиг меньших значений функции при меньшем числе вызовов оракулов.


Отметим также, что если функция $F(x,y)$~--- $\mu$-сильно выпуклая по совокупности переменных, то функция $g(y) = \min_{x\in Q_x} F(x,y)$ также будет $\mu$-сильно выпуклая. Более того, все это можно сформулировать в терминах $(\delta,\mu,L)$-оракула (см. \cite{gasnikov2018book} и цитированную там литературу). При $\mu = 0$ это сделано в работе \cite{Gas2015Equi}, при $\mu > 0$ доказательство практически дословно повторяет утверждения 1, 3 из \cite{Gas2015Equi} (см. также \cite{Gasnikov20Models}). Приведенное наблюдение позволяет обоснованно (с теоретической проработкой) использовать для решения внутренней задачи метод Вайды, а для решения внешний задачи использовать, например, быстрый градиентный метод. Однако такой подход будет предпочтительнее рассмотренного в данной статье только при весьма специальных (как правило, трудно выполнимых) условиях \cite{Gladin2020}.

\appendix{1}  

\begin{proofofstatement}{\ref{max_dot_product}}
Рассмотрим произвольный $y \in Q_y$.
\begin{equation}\label{dot_product}
    \left\langle \nabla g(\tilde{y}), \tilde{y}-y\right\rangle = \left\langle \nabla g(\tilde{y}) - \nabla g(y_*), \tilde{y}-y\right\rangle + \left\langle \nabla g(y_*), \tilde{y}-y\right\rangle
\end{equation}
Оценим сверху первое слагаемое, используя неравенство Коши--Буняковского и определение липшицевости градиента:
\begin{equation}\label{first_term}
    \left\langle \nabla g(\tilde{y}) - \nabla g(y_*), \tilde{y}-y\right\rangle \leq \left\| \nabla g(\tilde{y}) - \nabla g(y_*) \right\|_2 \left\| \tilde{y}-y \right\|_2 \leq L \left\| \tilde{y} - y_* \right\|_2 \left\| \tilde{y}-y \right\|_2.
\end{equation}
Из сильной выпуклости следует
\begin{equation*}
    g(\tilde{y}) \geq g(y_*)+\langle\nabla g(y_*), \tilde{y}-y_*\rangle+\frac{\mu}{2} \|\tilde{y}-y_*\|_2^{2}
\end{equation*}
Воспользовавшись неравенствами $g(\tilde{y}) - g(y_*) \leq \varepsilon$ и $\langle\nabla g(y_*), y-y_*\rangle \geq 0\ \forall y \in Q_y$, получим
\begin{equation}\label{argument_residue}
    \|\tilde{y}-y_*\|_2 \leq \sqrt{\frac{2 \varepsilon}{\mu}} \stackrel{\eqref{first_term}}{\Longrightarrow} \left\langle \nabla g(\tilde{y}) - \nabla g(y_*), \tilde{y}-y\right\rangle \leq L \left\| \tilde{y}-y \right\|_2 \sqrt{\frac{2 \varepsilon}{\mu}}.
\end{equation}
Теперь оценим сверху второе слагаемое в \eqref{dot_product}.
\begin{equation*}
    \left\langle \nabla g(y_*), \tilde{y}-y\right\rangle = \left\langle \nabla g(y_*), \tilde{y}-y_*\right\rangle + \left\langle \nabla g(y_*), y_*-y\right\rangle
\end{equation*}
Снова воспользовавшись критерием оптимальности точки $y_*$ и неравенством Коши-Буняковского, получим
\begin{equation*}
    \left\langle \nabla g(y_*), \tilde{y}-y\right\rangle \leq \| \nabla g(y_*) \|_2 \|\tilde{y}-y_*\|_2 \stackrel{\eqref{argument_residue}}{\leq} \| \nabla g(y_*) \|_2 \sqrt{\frac{2 \varepsilon}{\mu}}.
\end{equation*}
Объединив верхние оценки для обоих слагаемых, получим
\begin{equation*}
    \left\langle \nabla g(\tilde{y}), \tilde{y}-y\right\rangle \leq \left(L \left\| \tilde{y}-y \right\|_2 + \| \nabla g(y_*) \|_2 \right) \sqrt{\frac{2 \varepsilon}{\mu}},
\end{equation*}
откуда следует доказываемое утверждение.
\end{proofofstatement}

\begin{proofoftheorem}{\ref{delta_subgradient}}
Зафиксировав $x \in Q_x$, применим утверждение~\ref{max_dot_product} к функции $g(y):= F(x, y)$ и утверждение \ref{statement_gas}.
\end{proofoftheorem}

\begin{proofoftheorem}{\ref{theorem_approach_2}}
Согласно \eqref{estimate_restat_alg},
алгоритм \ref{alg:restart} сходится линейно, поэтому можно считать, что вспомогательная задача $\min_{y \in Q_y} F(x, y)$ решается сколь угодно точно за время $\widetilde{\mathcal{O}}\left(\sqrt{\frac{L}{\mu}} \right)$. Согласно теореме \ref{delta_subgradient}, это позволяет использовать $\delta$-субградиент, где $\delta$ убывает со скоростью геометрической прогрессии. Для внешней задачи используется метод Вайды, который также сходится линейно и имеет сложность $\widetilde{\mathcal{O}}\left(d \right)$. Таким образом, для решения задачи \eqref{main4} достаточно $\widetilde{\mathcal{O}} \left( d \right)$ вычислений $\partial_x F$ и обращений матриц размера $d \times d$, а также $\widetilde{\mathcal{O}} \left(d \sqrt{\frac{L}{\mu}} \right)$ вычислений $\nabla_y F$.
\end{proofoftheorem}

\begin{proofoftheorem}{\ref{theorem_approach_3}}
Согласно теореме \ref{theorem_estimate_varag}, Varag сходится линейно, поэтому можно считать, что вспомогательная задача $\min_{y \in Q_y} F(x, y)$ решается сколь угодно точно за время $\widetilde{\mathcal{O}}\left(m+\sqrt{\frac{m L}{\mu}} \right)$. Согласно теореме \ref{delta_subgradient}, это позволяет использовать $\delta$-субградиент, где $\delta$ убывает со скоростью геометрической прогрессии. Для внешней задачи используется метод Вайды, который также сходится линейно и имеет сложность $\widetilde{\mathcal{O}}\left(d \right)$ итераций. На каждой его итерации необходимо вычислять субградиенты всех $m$ слагаемых $\partial_x F_i$. Таким образом, для решения задачи достаточно $\widetilde{\mathcal{O}}\left(md\right)$ вычислений $\partial_x F_i$, $\widetilde{\mathcal{O}}\left(d\right)$ обращений матриц размера $d \times d$ и $\widetilde{\mathcal{O}} \left(dm+d\sqrt{\frac{m L}{\mu}} \right)$ вычислений $\nabla_y F_i$.
\end{proofoftheorem}

\AdditionalInformation{Гладин Е.Л.}{Московский физико-технический институт, Сколковский институт науки и технологий, студент магистратуры, Долгопрудный}{gladin.el@phystech.edu}

\AdditionalInformation{Алкуса М.}{Московский физико-технический институт, м.н.с., Долгопрудный}{mohammad.alkousa@phystech.edu}

\AdditionalInformation{Гасников А.В.}{Московский физико-технический институт, профессор, Долгопрудный}{gasnikov.av@mipt.ru}




\begin{thebibliography}{10}







\bibitem{Gasnikov20Models}
{\it Гасников А.В., Гасникова Е.В.}
Модели равновесного распределения транспортных потоков в больших сетях: учебное пособие. Москва: МФТИ, 2020.


\bibitem{vaidya1989}
{\it Vaidya P. M.}
A new algorithm for minimizing convex functions over convex
sets // In Foundations of Computer Science, 1989., 30th Annual Symposium, 1989. P. 338--343.


\bibitem{vaidya1989new}
{\it Vaidya P. M.}
A new algorithm for minimizing convex functions over convex sets // Mathematical Programming 73, Springer, 1996. P. 291--341.


\bibitem{nesterov2018lectures}
{\it Nesterov Yu.}
Lectures on convex optimization. Switzerland: Springer Optimization and Its Applications, 2018.
	
\bibitem{gasnikov2018book}
{\it Гасников А.В.} Современные численные методы оптимизации. Метод универсального градиентного спуска М.: МЦНМО, 2020.  

	
\bibitem{Nesterov1983}	
{\it Нестеров Ю.Е.}
Метод минимизации выпуклых функций со скоростью сходимости $O(1/k^2)$ // Докл. АН СССР. 1983. Т. 269. № 3. С. 543--547.



\bibitem{Lan_book2020}
{\it Lan G.}
First-order and Stochastic Optimization Methods for Machine Learning. Atlanta: Springer, 2020.


\bibitem{Lan_NIPSpaper}
{\it Lan G., Zhize Li, Yi Zhou}
A unified variance-reduced accelerated gradient
method for convex optimization. 33rd Conference on Neural Information Processing Systems (NeurIPS 2019), Vancouver, Canada. \url{https://arxiv.org/pdf/1905.12412.pdf}


\bibitem{Holderian2020}
{\it Bolte J., Glaudin L., Pauwels E., Serrurier M.}
A H\"{o}lderian backtracking method for min-max and min-min problems. \url{https://arxiv.org/pdf/2007.08810.pdf}



\bibitem{Alkousa2020}
{\it M. S. Alkousa, D. M. Dvinskikh, F. S. Stonyakin, A. V. Gasnikov, D. Kovalev}
Accelerated Methods for Saddle Point Problems// Comput. Math. and Math. Phys., 2020, Vol. 60, No. 11, P. 1787--1809.


\bibitem{Gladin2020}
{\it  Gladin E.,  Kuruzov I.,  Stonyakin F., Pasechnyuk D.,
Alkousa M.,  Gasnikov A.}
Solving strongly convex-concave composite saddle
point problems with a small dimension of one of the variables. \url{https://arxiv.org/pdf/2010.02280.pdf}



\bibitem{Jordan2020}
{\it  Tianyi L., Chi J., Michael. I. J.}
Near-Optimal Algorithms for Minimax Optimization. \url{https://arxiv.org/pdf/2002.02417v5.pdf}



\bibitem{Yuanhao2020}
{\it  Yuanhao W., Jian L.}
Improved Algorithms for Convex-Concave Minimax Optimization. \url{https://arxiv.org/pdf/2006.06359.pdf}


\bibitem{Zhongruo2020}
{\it  Zhongruo Wang, Krishnakumar Balasubramanian, Shiqian Ma, Meisam Razaviyayn}
Zeroth-Order Algorithms for Nonconvex Minimax Problems with Improved Complexities. \url{https://arxiv.org/pdf/2001.07819.pdf}



\bibitem{Jungers2011}
{\it  Jungers M.,  Tr\'{e}lat E.,  Abou-Kandil H.}
Min-max and min-min Stackelberg strategies with closed-loop information structure // Journal of Dynamical and Control Systems, Springer Verlag, 2011, 17 (3), P.387--425. 



\bibitem{Konur2017}
{\it   Konur D., Farhangi H.}
Set-based Min-max and Min-min Robustness for Multi-objective
Robust Optimization // Proceedings of the 2017 Industrial and Systems Engineering Research Conference K. Coperich, E. Cudney, H. Nembhard, eds. 

\bibitem{Bubeck}
{\it   Bubeck S.}
Convex optimization: algorithms and complexity // Foundations and Trends in Machine Learning. 2015. Vol. 8, №3–4. P. 231--357.


\bibitem{FGM_Gasnikov2019}
{\it   Tyurin A. I., Gasnikov A.V.}
Fast gradient descent method for convex optimization problems with
an oracle that generates a $(\delta, L)$-model of a function in a requested point // Comput. Math. and Math. Phys., 2019. Vol. 59, No. 7, P. 1137–1150.

\bibitem{Gas2015Equi}
{\it Гасников А.В., Двуреченский П.Е., Камзолов Д.И., Нестеров Ю.Е., Спокойный В.Г., Стецюк П.И., Суворикова А.Л., Чернов А.В.}
Поиск равновесий в многостадийных транспортных моделях // Труды Московского физико-технического института, 7.4(28) (2015).

\bibitem{polyak1983intro}
{\it Поляк Б.Т.}
Введение в оптимизацию // М.: Наука, 1983. 384 с.






\end{thebibliography}
\end{document}